\documentclass[11pt]{article}
\usepackage{amscd}
\usepackage{amsfonts}
\usepackage{amsmath}
\usepackage{amssymb}
\usepackage{amsthm}
\usepackage{bbm}
\usepackage{CJK}
\usepackage{fancyhdr}
\usepackage{graphicx}
\usepackage{hyperref}
\usepackage{indentfirst}
\usepackage{latexsym}
\usepackage{mathrsfs}
\usepackage{xypic}

\usepackage[top=1in,bottom=1in,left=1.25in,right=1.25in]{geometry}
\newtheorem{theorem}{Theorem}[section]
\newtheorem{lemma}[theorem]{Lemma}
\newtheorem{definition}[theorem]{Definition}
\newtheorem{proposition}[theorem]{Proposition}
\newtheorem{example}[theorem]{Example}

\newtheorem{corollary}[theorem]{Corollary}

\def\<{\langle}
\def\>{\rangle}
\def\a{\alpha}
\def\b{\beta}

\def\c{\cdot}

\def\g{\gamma}

\def\o{\otimes}

\date{}
\begin{document}
\renewcommand{\baselinestretch}{1.2}
\renewcommand{\arraystretch}{1.0}
\title{\bf The constructions of   3-Hom-Lie bialgebras }
\date{}
\author{{\bf Shuangjian Guo$^{1}$, Xiaohui Zhang$^{2}$,  Shengxiang Wang$^{3}$\footnote
        { Corresponding author(Shengxiang Wang):~~wangsx-math@163.com} }\\
{\small 1. School of Mathematics and Statistics, Guizhou University of Finance and Economics} \\
{\small  Guiyang  550025, P. R. of China} \\
{\small 2.  School of Mathematical Sciences, Qufu Normal University}\\
{\small Qufu  273165, P. R. of China}\\
{\small 3.~ School of Mathematics and Finance, Chuzhou University}\\
 {\small   Chuzhou 239000,  P. R. of China}}
 \maketitle
\begin{center}
\begin{minipage}{12.cm}
\begin{center}{\bf ABSTRACT}\end{center}
In this paper, we first introduce the notion  of a 3-Hom-Lie bialgebra and prove that it is equivalent to a  Manin triple of 3-Hom-Lie algebras.
 Also, we study the $\mathcal{O}$-operator and construct solutions of the  3-Lie classical Hom-Yang-Baxter equation interms of $\mathcal{O}$-operators and 3-Hom-pre-Lie algebras.    Finally, we  show that a 3-Hom-Lie algebra has a
phase space if and only if it is sub-adjacent to a 3-Hom-pre-Lie algebra.
\medskip

\medskip

\end{minipage}
\end{center}

\section{ Introduction}

3-Lie algebras are special types of $n$-Lie algebras and have close relationships with many important fields in mathematics
and mathematical physics \cite{BL08}-\cite{BL8}. The structure of 3-Lie algebras is closely linked to the supersymmetry and gauge
symmetry transformations of the world-volume theory of multiple coincident $M_2$-branes and is applied to the study of
the Bagger-Lambert theory. Moreover, the $n$-Jacobi identity can be regarded as a generalized Plucker relation in the physics
literature.  In particular, the metric 3-Lie algebras, or more generally, the 3-Lie algebras with
invariant symmetric bilinear forms attract even more attentions in physics. Recently, many
more properties and structures of 3-Lie algebras have been developed, see  [5-7], [9], [16]  and references cited therein.

The notion of Hom-Lie algebras was introduced by Hartwig, Larsson, and Silvestrov in \cite{HL06} as part of
a study of deformations of the Witt and the Virasoro algebras. In a Hom-Lie algebra, the Jacobi identity
is twisted by a linear map, called the Hom-Jacobi identity. Some $q$-deformations of the Witt and the
Virasoro algebras have the structure of a Hom-Lie algebra \cite{HL06}. Because of close relation to discrete and
deformed vector fields and differential calculus \cite{GZW}-\cite{H99}. More people
pay special attention to this algebraic structure. In particular, representations and cohomologies
of Hom-Lie algebras were studied in \cite{S12},   the cohomologies adapted to central extensions and
deformations are studied in \cite{AMM11}, construction of 3-Hom-Lie algebras from Hom-Lie algebras are studied in \cite{AMS11},
construction of  Rota-Baxter multiplicative 3-ary Hom-Nambu-Lie algebras  in \cite{SC15} and extensions of 3-Hom-Lie algebras are studied in \cite{LC15},  construction of 3-Lie classical Hom-Yang-Baxter equation in \cite{WW17}. This provides a good starting point for discussion and further research.

Looking forward, further attempts conld prove quite beneficial to the literature. Motivated by these results of \cite{BC16}, \cite{SB14} and \cite{SC13}, the paper is organized as follows.
 In Section 2,  we recall several concepts and results, and  introduce the notions of the matched pair of 3-Hom-Lie algebras, the 3-Hom-Lie bialgebra and the Manin triple of 3-Hom-Lie algebras.
  The following three expressions are equivalent: $(L, L^{\ast})$ is a 3-Hom-Lie bialgebra, $(L, L^{\ast}, ad^{\ast}, a\partial^{\ast}, \a, \a^{\ast})$ is a matched pair of 3-Hom-Lie algebras, and  $((L\oplus L^{\ast}, (\c,\c),\a+\a^{\ast}) L, L^{\ast})$  is a standard Manin triple of 3-Hom-Lie algebras.
   In Section 3,  we introduce the notion of the $\mathcal{O}$-operator and construct solutions of  the  3-Lie classical Hom-Yang-Baxter equation in terms of $\mathcal{O}$-operators and 3-Hom-pre-Lie algebras.
  In Section 4, we study representations of 3-Hom-pre-Lie algebras. In Section 5,
we introduce the notion of the phase space of a 3-Hom-Lie algebra and show that a 3-Hom-Lie algebra has a
phase space if and only if it is sub-adjacent to a 3-Hom-pre-Lie algebra.

\section{ Matched pairs, Manin triples and  3-Hom-Lie bialgebras }
\def\theequation{\arabic{section}.\arabic{equation}}
\setcounter{equation} {0}

In this section, we recall some basic notions and facts about 3-Hom-Lie algebra  and present some examples.
Then we give an equivalent description of 3-Hom-Lie bialgebras,  matched pair and Manin triples of 3-Hom-Lie algebras.
\medskip

\begin{definition} (\cite{AMM11})
A 3-Hom-Lie algebra is a triple $(L, [\c, \c, \c], \a)$ consisting of a vector space $L$, a 3-ary skew-symmetric
operation $[\c,\c,\c]_L: \wedge^{3}L\rightarrow L$  and a linear map $\a:L\rightarrow L$ satisfying the following Hom-Jacobi identity
\begin{eqnarray*}
[\a(x), \a(y), [u,v,w]]&=&[[x,y,u],\a(v),\a(w)] +[\a(u), [x,y,v], \a(w)]\nonumber \\
&&+[\a(u), \a(v),[x,y,w]],
\end{eqnarray*}
for any $x, y, u, v,w \in L$.
\end{definition}

 A 3-Hom-Lie algebra is called a multiplicative 3-Hom-Lie algebra if $\alpha$ is an algebra morphism, i.e. for any $x, y, z \in  L$,
$\alpha([x, y, z]) = [\alpha(x), \alpha(y), \alpha(z)]$. A 3-Hom-Lie algebra is called regular if $\alpha$ is an algebra automorphism.

\begin{example}
Let $(L, [\c,\c,\c])$ be a 3-Lie algebra and $\a: L\rightarrow L$ an algebra morphism,
 then the algebra $(L, [\c,\c,\c]_{\a}, \a)$ is a 3-Hom-Lie algebra, where $[\c,\c,\c]_{\a}$ is defined by
\begin{eqnarray*}
[x_1,x_2,x_3]_{\a}=[\a(x_1), \a(x_2), \b(x_3)].
\end{eqnarray*}
\end{example}

\begin{example}
Let $(L, [\c,\c,\c], \a)$ be a 3-Hom-Lie algebra and $\b: L\rightarrow L$ an algebra morphism such that $\a\b=\b\a$,
 then  $(L, [\c,\c,\c]_{\a\b}=[\c,\c,\c]\circ (\b\o \b\o \a), \a\circ \b)$ is a 3-Hom-Lie algebra.
\end{example}

\begin{example}
Let $(L, [\c,\c,\c], \a)$ be a 3-Hom-Lie algebra over a filed $\mathbb{F}$  and $t$ an indeterminate,
 define $\hat{L}: \{\sum (x\o t+y\o t^n)\subset L\o (\mathbb{F}[t]/t^{n+1})|x,y\in L)\}$, $\hat{\a}(\hat{L})=\{\sum (\a(x)\o t+\a(y)\o t^n):x,y\in L\}.$
 Then $(\hat{L}, \hat{\a})$  ia a 3-Hom-Lie algebra with the operation $[x_1\o t^{i_1}, x_2\o t^{i_2},x_3\o t^{i_3}]=[x_1,x_2,x_3]\o t^{\sum i_j}$ for all $x_1,x_2,x_3\in L$ and  $i_1,i_2,i_3\in \{1,2,3\}$.

\end{example}

\begin{definition} (\cite{LC15})
A representation of a multiplicative 3-Hom-Lie algebra $(L, [\c, \c, \c], \a)$ on the vector space $V$ with respect to
$A\in gl(V)$ is a bilinear map $\rho: L \wedge L \rightarrow gl(V)$, such that for any $x, y, z, u, \in L$, the following equalities are satisfied
\begin{eqnarray*}
&&\rho(\a(u), \a(v))\circ A=A\circ \rho(u, v),\\
&& \rho([x,y,z], \a(u))\circ A=\rho(\a(y), \a(z))\rho(x,u)+\rho(\a(z), \a(x))\rho(y,u)\\
&&\hspace{4cm}+\rho(\a(x), \a(y))\rho(z,u),\\
&& \rho(\a(x), \a(y))\rho(z,u)=\rho(\a(z), \a(u))\rho(x,y)+\rho([x,y,z], \a(u))\circ A\\
&&\hspace{4cm}+\rho( \a(z), [x,y,u],)\circ A.
\end{eqnarray*}
Then $(V, \rho, A)$ is called a representation of $L$, or $V$ is an $L$-module.
\end{definition}

\begin{lemma}(\cite{LC15})
Let $(V, \rho, A)$ be a representation of a 3-Hom-Lie algebra $(L, [\c,\c,\c], \a)$.
Then $(V^{\ast}, \rho^{\ast}, A^{\ast})$ is also a representation of the 3-Hom-Lie algebra $(L, [\c,\c,\c], \a)$,
 which is called the dual representation.
\end{lemma}

\begin{lemma}(\cite{LC15})
 Let $(V, \rho, A)$ be a representation of a 3-Hom-Lie algebra $(L, [\c,\c,\c], \a)$.
 Then there is a 3-Hom-Lie algebra structure on the direct sum of vector spaces $L\oplus V$, defined by
\begin{eqnarray*}
[x_1+v_1, x_2+v_2, x_3+v_3]=[x_1,x_2,x_3] +\rho(x_1,x_2)v_3+\rho(x_2,x_3)v_1+\rho(x_3,x_`)v_2,\\
(\a+A)(x+v)=\a(x)+Av,
\end{eqnarray*}
for any $x_1,x_2,x_3,x\in L$ and $v_1,v_2,v_3,v\in V$.
\end{lemma}

\begin{definition} (\cite{AMM11})
A linear map $D : L \rightarrow  L$  is called a derivation of the multiplicative
3-Hom-Lie algebra $(L, [\c,\c,\c], \a)$, if
\begin{eqnarray*}
D\circ \a=\a\circ D,~~
D[x,y,z]=[D(x), y, z]+[x, D(y), z]+[x, y, D(z)],
\end{eqnarray*}
for any $x,y,z\in L.$ Denoted by $Der(L)$.
\end{definition}

\begin{proposition}
Let $(L, [\c,\c,\c], \a)$ and $(L', [\c,\c,\c]', \a')$ be 3-Hom-Lie algebras. Suppose that there are skew-symmetric linear maps $\rho: L\o L\rightarrow gl(L')$ and $\mu: L'\o L'\rightarrow gl(L)$  which are representations of $L$ and $L'$ respectively,
 satisfying the following equations:
\begin{eqnarray}
&& \mu(\a'(a_4), \a'(a_5))[x_1,x_2,x_3]-[\mu(a_4,a_5)x_1, \a(x_2),\a(x_3)]\nonumber\\
&&\hspace{2cm}-[\a(x_1), \mu(a_4,a_5)x_2,\a(x_3)]-[\a(x_1),\a(x_2), \mu(a_4,a_5)x_3]=0,
\end{eqnarray}
\begin{eqnarray}
&& \mu(\rho(x_1,x_4)a_5,a_3)\a(x_2)-\mu(\rho(x_2,x_4)a_5,a_3)\a(x_1)\nonumber\\
&&\hspace{2cm}-\mu(\rho(x_1,x_2)a_3,a_5)\a(x_4)+[\a(x_1),\a(x_2), \mu(a_3,a_5)x_4]=0,\\
&& [\mu(a_2,a_3)x_1,\a(x_4), \a(x_5)]- \mu(\a'(a_2), \a'(a_3))[x_1,x_4,x_5]\nonumber\\
&&\hspace{2cm}-\mu(\rho(x_4,x_5)a_2,a_3)\a(x_1)-\mu(a_2,\rho(x_4,x_5)a_3)\a(x_1)=0,\\
&& \rho(\a(x_4), \a(x_5))[a_1,a_2,a_3]'-[\rho(x_4,x_5)a_1, \a'(a_2),\a'(a_3)]'\nonumber\\
&&\hspace{2cm}-[\a'(a_1), \rho(x_4,x_5)a_2,\a'(a_3)]'-[\a'(a_1),\a'(a_2), \rho(x_4,x_5)a_3]'=0,\\
&& \rho(\mu(a_1,a_4)x_5,x_3)\a'(a_2)-\rho(\mu(a_2,a_4)x_5,x_3)\a'(a_1)\nonumber\\
&&\hspace{2cm}-\rho(\mu(a_1,a_2)x_3,x_5)\a'(a_4)+[\a'(a_1),\a'(a_2), \rho(x_3,x_5)a_4]'=0,\\
&& [\rho(x_2,x_3)a_1,\a'(a_4), \a'(a_5)]- \rho(\a(x_2), \a(x_3))[a_1,a_4,a_5]'\nonumber\\
&&\hspace{2cm}-\rho(\mu(a_4,a_5)x_2,x_3)\a'(a_1)-\rho(x_2,\mu(a_4,a_5)x_3)\a'(a_1)=0,
\end{eqnarray}
for any $x_i\in L$ and $a_i\in L', 1\leq i\leq 5$. Then there is a 3-Hom-Lie algebra structure on $L\oplus L'$ defined by
\begin{eqnarray*}
&& (\a+\a')(x_1+a_1)=\a(x_1)+\a'(a_1),\\
&&[x_1+a_1,x_2+a_2,x_3+a_3]_{L\oplus L'}= [x_1,x_2,x_3]+\rho(x_1,x_2)a_3+\rho(x_3,x_1)a_2+\rho(x_2,x_3)a_1\\
&&\hspace{5cm}= [a_1,a_2,a_3]'+\mu(a_1,a_2)x_3+\mu(a_3,a_1)x_2+\mu(a_2,a_3)x_1,
\end{eqnarray*}
for any $x_i\in L$ and $a_i\in L', 1\leq i\leq 5$.
Moreover, $(L, [\c,\c,\c],  L', [\c,\c,\c]', \rho, \mu,  \a, \a')$ satisfying the above conditions is called a matched pair of 3-Hom-Lie algebras.
\end{proposition}
{\bf Proof.} Similar to \cite{DB18}.  \hfill $\square$

\begin{definition}
Let $(L, [\c, \c, \c], \a)$ be a 3-Hom-Lie algebra. A bilinear form $(\c, \c)_A$ on $L$ is called invariant
if it satisfies
\begin{eqnarray}
([x,y,z], \a(u))_A+([x,y,u], \a(z))_A=0,~~~~\forall x,y,z,u\in L.
\end{eqnarray}
A 3-Hom-Lie algebra $L$ is called pseudo-metric if there is a non-degenerate symmetric invariant
bilinear form on $L$.
\end{definition}

\begin{definition}
A Manin triple of 3-Hom-Lie algebras consists of a pseudo-metric 3-Hom-Lie
algebra $(L, (\c, \c)_A, \a)$ and 3-Hom-Lie algebras $L_1$ and $L_2$ such that

(1) $L_1, L_2$ are isotropic 3-Hom-Lie subalgebras of $L$;

(2) $L = L_1 \oplus L_2$ as the direct sum of vector spaces;

(3) For all $x_1, y_1 \in  L_1$ and $x_2, y_2 \in  L_2$, we have pr$_1
[x_1, y_1, x_2] = 0$ and pr$_2[x_2, y_2, x_1] =0$,
where pr$_1$ and pr$_2$ denote the projections from $L_1 \oplus L_2$ to $L_1, L_2$, respectively.
\end{definition}

Let $(L, [\c, \c, \c]_L, \a)$  be a 3-Hom-Lie algebra and $(L^\ast, [\c, \c, \c]^{\ast}, \a^\ast)$ its dual.
There is a natural nondegenerate symmetric bilinear form $(\c,\c)$ on $L \oplus L^{\ast}$  given by
\begin{eqnarray}
(x+\xi,y+\eta)=<x,\eta>+<\xi, y>,~~\forall x,y\in L, \xi, \eta\in L^{\ast}.
\end{eqnarray}
There is also a bracket operation $[\c, \c, \c]_{L\oplus L^{\ast}}$ on $L \oplus L^{\ast}$ given by
\begin{eqnarray}
&&(\a+\a^{\ast})(x+\xi)=\a(x)+\a^{\ast}(\xi),\nonumber\\
&& [x+\xi,y+\eta, z+\g]_{L\oplus L^{\ast}}=[x,y,z]+ad^{\ast}_{x,y}\g +ad^{\ast}_{y,z}\xi+ad^{\ast}_{z,x}\eta \nonumber\\
&&\hspace{4cm}+ a\partial^{\ast}_{\xi,\eta}z +a\partial^{\ast}_{\eta, \gamma}x +a\partial^{\ast}_{\gamma,\xi}y+[\xi, \eta, \g]^{\ast},
\end{eqnarray}
where $ad^{\ast}$ and $a\partial^{\ast}$ are the coadjoint representations of $L$ on $L^{\ast}$ and $L^{\ast}$ on $L$, respectively.
Note that the bracket operation $[\c, \c, \c]_{L\oplus L^{\ast}}$ is naturally invariant with respect to the
symmetric bilinear form $(\c,\c)$, and satisfies the condition (3) in Definition 2.11.
If $(L\oplus L^{\ast}, [\c, \c, \c]_{L\oplus L^{\ast}}, \a+\a^{\ast})$  is a 3-Hom-Lie algebra, then obviously $L$ and $L^{\ast}$ are isotropic subalgebras.
Consequently, $((L\oplus L^{\ast}, (\c,\c),\a+\a^{\ast}),  L, L^{\ast})$  is a Manin triple, which is called the standard Manin triple of 3-Hom-Lie algebras.
\medskip

We immediately get the following relation between matched pairs and Manin triple of 3-Hom-Lie algebras.

\begin{proposition}
Let $(L, [\c, \c, \c]_L, \a)$  be a 3-Hom-Lie algebra and $(L^\ast, [\c, \c, \c]^{\ast}, \a^\ast)$ its dual.
Then $((L\oplus L^{\ast}, (\c,\c),\a+\a^{\ast}) L, L^{\ast})$  is a standard Manin triple if and only if $(L, L^{\ast}, ad^{\ast},\\
 a\partial^{\ast}, \a, \a^{\ast})$ is a matched pair.
\end{proposition}

\begin{lemma}
Let $(L, [\c, \c, \c]_L, \a)$  be a 3-Hom-Lie algebra and $(L^\ast, [\c, \c, \c]^{\ast}, \a^\ast)$ its dual.
 Then $(L, L^{\ast}, ad^{\ast},a\partial^{\ast}, \a, \a^{\ast})$ is a matched pair if and only if (2.1), (2.2) and (2.3) hold for $\rho=ad^{\ast}$ and $\mu=a\partial^{\ast}$.
\end{lemma}
{\bf Proof.} $\Leftarrow$, For any $x_1,x_2,x_4\in L$ and $a_3,a_5,a_6\in L^{\ast}$, we have
\begin{eqnarray*}
&&<-a\partial^{\ast}_{ad^{\ast}_{x_1,x_2}a_3,a_5}\a(x_4)+a\partial^{\ast}_{ad^{\ast}_{x_1,x_4}a_5,a_3}\a(x_2)-
a\partial^{\ast}_{ad^{\ast}_{x_2,x_4}a_5,a_3}\a(x_1)\\
&&+[\a(x_1),\a(x_2), a\partial^{\ast}_{a_3,a_5}x_4], a_6>\\
&=& <[ad^{\ast}_{x_1,x_2}a_3,a_5, a_6]^{\ast}, \a(x_4)>-<[ad^{\ast}_{x_1,x_4}a_5,a_3, a_6]^{\ast}, \a(x_2)>\\
&& +<[ad^{\ast}_{x_2,x_4}a_5,a_3, a_6]^{\ast}-ad^{\ast}_{\a(x_2),a\partial^{\ast}_{a_3,a_5}x_4 }a_6, \a(x_1)>\\
&=& <x_1, ad^{\ast}_{x_2,a\partial^{\ast}_{a_5,a_6}x_4 }\a^{\ast}(a_3)+ad^{\ast}_{x_4,a\partial^{\ast}_{a_3,a_6}x_2 }\a^{\ast}(a_5)\\
&&+ad^{\ast}_{x_2,x_4}a_5,\a^{\ast}(a_3), \a^{\ast}(a_6)]^{\ast}-ad^{\ast}_{x_2,a\partial^{\ast}_{a_3,a_5}x_4 }\a^{\ast}(a_6)>,
\end{eqnarray*}
which implies the equivalence between (2.2) and (2.5). The proofs of the others cases are similar.
\begin{theorem}
Let $(L, [\c, \c, \c]_L, \a)$ be a 3-Hom-Lie algebra and $\Delta:A\rightarrow A\o A\o A$ a linear map. Suppose that $\Delta^{\ast}: A^{\ast}\o A^{\ast}\o A^{\ast}\rightarrow A^{\ast}$ defines a 3-Hom-Lie algebra structure $[\c,\c,\c]^{\ast}$ on $L^{\ast}$. Then $(L, L^{\ast}, ad^{\ast}, a\partial^{\ast}, \a, \a^{\ast})$ is a matched pair if and only if  the following equations are satisfied:
\begin{eqnarray}
&& \Delta([x,y,z])=(\a\o\a\o ad_{y,z}) \Delta(x)+(\a\o\a\o ad_{z,x}) \Delta(y)+(\a\o\a\o ad_{x,y}) \Delta(z),~~~\\
&& \Delta([x,y,z])=(\a\o\a\o ad_{y,z}) \Delta(x)+(\a\o ad_{y,z} \o \a) \Delta(x)+(ad_{y,z}\o \a\o\a) \Delta(x),~~~\\
&& (ad_{x,y}\o \a\o\a+\a\o\a\o ad_{x,y})\Delta(z)=(\a\o ad_{z,x} \o \a) \Delta(y)+(\a\o ad_{y,z} \o \a) \Delta(x),~~~~~~
\end{eqnarray}
for any $x,y,z\in L$.
\end{theorem}
{\bf Proof.} Let $\{e_1,e_2,...,e_n\}$ be a basis of $L$ and $\{e^{\ast}_1,e^{\ast}_2,...,e^{\ast}_n\}$ the dual basis. Suppose
\begin{eqnarray*}
[e_i,e_j,e_k]=\sum_{l=1}^{n}c_{ijk}^{l}e_l,    [e^{\ast}_i,e^{\ast}_j,e^{\ast}_k]^{\ast}=\sum_{l=1}^{n}d_{ijk}^{l}e^{\ast}_l.
\end{eqnarray*}
Let
\begin{eqnarray*}
&&\a(e_i)=\sum_{s}f_se_s, \a(e_j)=\sum_{n}g_ne_n,  \a(e_k)=\sum_{n}h_me_m,\\
&& \a^{\ast}(e^{\ast}_{\xi})=\sum_{s}f^{\ast}_se^{\ast}_s, \a^{\ast}(e^{\ast}_\eta)=\sum_{n}g^{\ast}_ne^{\ast}_n,  \a^{\ast}(e^{\ast}_k)=\sum_{m}h^{\ast}_me^{\ast}_m.
\end{eqnarray*}
Then we have
\begin{eqnarray*}
ad^{\ast}_{e_i,e_j}e_k^{\ast}=-\sum_{l=1}^{n}c_{ijk}^{l}e_l^{\ast}, a\partial^{\ast}_{e^{\ast}_i,e^{\ast}_j}e_k=-\sum_{l=1}^{n}d_{ijk}^{l}e_l,\Delta(e_k)=\sum^{n}_{i,j,l=1}d_{ijl}^{k}e_i\o e_j\o e_k.
\end{eqnarray*}
By (2.1), we have
\begin{eqnarray*}
&&a\partial^{\ast}_{\a^{\ast}(e^{\ast}_{\xi}), \a^{\ast}(e^{\ast}_{\eta})}[e_i,e_j,e_k]-[a\partial^{\ast}_{e^{\ast}_{\xi}, e^{\ast}_{\eta}}e_i,\a(e_j),\a(e_k)]\\
&&~~~~~~~-[\a(e_i), a\partial^{\ast}_{e^{\ast}_{\xi}, e^{\ast}_{\eta}}e_j,\a(e_k)]
-[\a(e_i), \a(e_j), a\partial^{\ast}_{e^{\ast}_{\xi}, e^{\ast}_{\eta}}e_k]=0,
\end{eqnarray*}
which gives
\begin{eqnarray*}
\sum_{l=1}^{n}(-f^{\ast}_sg^{\ast}_nd_{snm}^{l}c^{l}_{ijk}+g_nh_md^{i}_{\xi\eta l}c_{lnm}^{m}+f_sh_md^{j}_{\xi\eta l}c_{slm}^{m}+f_sg_nd^{k}_{\xi\eta l}c_{snl}^{m}=0,
\end{eqnarray*}
as the coefficient of $e_m$. On the other hand, the left hand side of the obove equation is also the coefficient of $e_{\xi}\o e_{\eta}\o e_m$ in (2.10).
Thus, we deduce that (2.1) is equivalent to (2.2).  The proofs of the others cases are similar.

\begin{definition}
Let $(L, [\c, \c, \c]_L, \a)$ be a 3-Hom-Lie algebra and $\Delta:A\rightarrow A\o A\o A$ a linear map. Suppose that $\Delta^{\ast}: A^{\ast}\o A^{\ast}\o A^{\ast}\rightarrow A^{\ast}$ defines a 3-Hom-Lie algebra structure $[\c,\c,\c]^{\ast}$ on $L^{\ast}$. If $\Delta$ satisfies (2.10) and (2.11), then we call $(L, \a, \Delta)$ a double construction 3-Hom-Lie bialgebra.
\end{definition}

Combining Proposition 2.9, Proposition 2.12 and   Theorem 2.14, we have

\begin{theorem}
Let $(L, [\c, \c, \c]_L, \a)$ be a 3-Hom-Lie algebra and $\Delta:A\rightarrow A\o A\o A$ a linear map. Suppose that $\Delta^{\ast}: A^{\ast}\o A^{\ast}\o A^{\ast}\rightarrow A^{\ast}$ defines a 3-Hom-Lie algebra structure $[\c,\c,\c]^{\ast}$ on $L^{\ast}$. Then the following statements are equivalent:

(1) $(L, \a, \Delta)$ is  a double construction 3-Hom-Lie bialgebra.

(2) $((L\oplus L^{\ast}, (\c,\c),\a+\a^{\ast}) L, L^{\ast})$  is a standard Manin triple of 3-Hom-Lie algebras.

(3) $(L, L^{\ast}, ad^{\ast}, a\partial^{\ast}, \a, \a^{\ast})$ is a matched pair of 3-Hom-Lie algebras.

\end{theorem}

\section{$\mathcal{O}$-operators and 3-Hom-pre-Lie algebras}
\def\theequation{\arabic{section}.\arabic{equation}}
\setcounter{equation} {0}

In this section, we mainly study the $\mathcal{O}$-operator of a 3-Hom-Lie algebra
 and present a kind of 3-Hom-Lie Yang-Baxter equation solutions.

\begin{definition}
Let $(L, [\c, \c, \c], \a)$ be 3-Hom-Lie algebra and $(V, \rho_A)$ a representation.
 A linear operator $T:V\rightarrow L$ is called an $\mathcal{O}$-operator associated to $(V, \rho_A)$ if $T$ satisfies: for any $u,v,w\in L$,
\begin{eqnarray}
&&  \a\circ T=T \circ A,\\
&& [Tu, Tv, Tw]_{L}=T(\rho_A(Tu, Tv)w+\rho_A(Tv, Tw)u+\rho_A(Tw, Tu)v).
\end{eqnarray}
\end{definition}

\begin{example}
Let $(L, [\c, \c, \c]_L, \a)$ be a 3-Hom-Lie algebra.
 An $\mathcal{O}$-operator of $L$ associated to the adjoint representation $(L,ad, \a)$
 is nothing but the Rota-Baxter operator of weight zero introduced by \cite{SC15}.
\end{example}

\begin{definition}
A 3-Hom-pre-Lie algebra  is a triple $(L, \{\c, \c, \c\}, \a)$ consisting of a vector space $L$,   with a linear map $\{\c, \c, \c\}: L\o L\o L\rightarrow L$ and a linear map $\a:L\rightarrow L$ satisfying
\begin{eqnarray}
&& \{x, y, z\}=-\{y, x, z\},\\
 \{\a(x), \a(y), \{z,u,w\}\}&=&\{[x,y,z]_L, \a(u), \a(v)\}+ \{\a(z), [x,y,u]_L, \a(v)\}\nonumber\\
&& +\{\a(z), \a(u), [x,y,v]_L\},\\
 \{[x,y,z]_L, \a(u), \a(v)\}&=&\{\a(x), \a(y), [z,u,v]_L\}+\{\a(y), \a(z), [x,u,v]_L\}\nonumber\\
&& +\{\a(z), \a(x), [y,u,v]_L\},
\end{eqnarray}
for any $x,y,z,u,v\in L$ and $[\c,\c,\c]_{C}$ is defined by
\begin{eqnarray}
[x,y,z]_{C}=\{x,y,z\}+\{y,z,x\}+\{z,x,y\}.
\end{eqnarray}
\end{definition}
\begin{proposition}
Let $(L, \{\c, \c, \c\}, \a)$ be a 3-Hom-pre-Lie algebra. Then the induced 3-commutator given by (3.6) defines a 3-Hom-Lie algebra.
\end{proposition}
{\bf Proof.} It is easy to check that $[\c,\c,\c]_C$ is skew-symmetric.  For any $x_1,x_2, x_3, x_4, x_5\in L$, we have
\begin{eqnarray*}
&&[\a(x_1), \a(x_2), [x_3, x_4, x_5]_C]_C-[[x_1, x_2, x_3]_{C},\a(x_4), \a(x_5)]_{C}-[\a(x_3), [x_1, x_2, x_4]_{C}, \a(x_5)]_{c}\\
&& -[\a(x_3), \a(x_4), [x_1, x_2, x_5]_{C},]_{C}\\
&=& \{\a(x_1), \a(x_2), \{x_3, x_4, x_5\}\}+\{\a(x_1), \a(x_2), \{x_4, x_5, x_3\}\}+\{\a(x_1), \a(x_2), \{x_5, x_3, x_4\}\}\\
&& +\{\a(x_2), [x_3, x_4, x_5]_{C}, \a(x_1)\}+\{ [x_3, x_4, x_5]_{L}, \a(x_1), \a(x_2)\}\\
&& -\{\a(x_4), \a(x_5), \{x_1, x_2, x_3\}\}-\{\a(x_4), \a(x_5), \{x_2, x_3, x_1\}\}-\{\a(x_4), \a(x_5), \{x_3, x_1, x_2\}\}\\
&& -\{[x_1, x_2, x_3]_{L},\a(x_4), \a(x_5)\}-\{\a(x_5), [x_1, x_2, x_3]_{C},\a(x_4)\}\\
&& -\{\a(x_5), \a(x_3), \{x_1, x_2, x_4\}\}-\{\a(x_5), \a(x_3), \{x_2, x_4, x_1\}\}-\{\a(x_5), \a(x_3), \{x_2, x_4, x_1\}\}\\
&& -\{[x_1, x_2, x_4]_{L},\a(x_5), \a(x_3)\}-\{\a(x_3), [x_1, x_2, x_4]_{C},\a(x_5)\}\\
&& -\{\a(x_3), \a(x_4), \{x_1, x_2, x_5\}\}-\{\a(x_3), \a(x_4), \{x_2, x_5, x_1\}\}-\{\a(x_3), \a(x_4), \{x_5, x_1, x_2\}\}\\
&& -\{[x_1, x_2, x_5]_{C},\a(x_3), \a(x_4)\}-\{\a(x_4), [x_1, x_2, x_5]_{C},\a(x_3)\}\\
&=& 0.
\end{eqnarray*}
Thus the proof is finished.\hfill $\square$

\begin{definition}
Let $(L, \{\c, \c, \c\}, \a)$ be a 3-Hom-pre-Lie algebra. The 3-Hom-Lie algebra $(L^c, [\c, \c, \c]_{C}, \a)$ is called the sub-adjacent 3-Hom-Lie algebra of $(L, \{\c, \c, \c\}, \a)$ and $(L, \{\c, \c, \c\}, \a)$ is called a compatible 3-Hom-pre-Lie algebra of the  3-Hom-Lie algebra $(L^c, [\c, \c, \c]_{C}, \a)$.
\end{definition}

Define the left multiplication $\mathcal{L} : \wedge^2 L\rightarrow gl(L)$ by $\mathcal{L} (x, y)z = \{x, y, z\}$  for all $x, y, z \in L$. Then
$(L,\mathcal{L}, \a)$ is a representation of the 3-Lie algebra $L$. Moreover, we define the right multiplication
$\mathcal{R} : \wedge^2 L\rightarrow gl(L)$ by $\mathcal{L} (x, y)z = \{z, x, y\}$. If there is a 3-pre-Lie algebra structure on its dual
space $L^{\ast}$, we denote the left multiplication and right multiplication by $\mathcal{L^{\ast}}$ and $\mathcal{R^{\ast}}$ respectively.

\begin{proposition}
Let $(L, [\c, \c, \c]_L, \a)$ be a 3-Hom-Lie algebra and $(V, \rho)$ a representation. Suppose that the linear map $T: V\rightarrow L$ is an $\mathcal{O}$-operator associated to $(V, \rho, \b)$. Then there exists a 3-Hom-pre-Lie algebra structure on $V$ given by
\begin{eqnarray*}
\{u,v,w\}=\rho(Tu,Tv)w, ~~~\forall u,v,w\in V.
\end{eqnarray*}
\end{proposition}
{\bf Proof.} For any $u,v,w\in V$,  we have
\begin{eqnarray*}
\{u,v,w\}=\rho(Tu,Tv)w=-\rho(Tv,Tu)w=-\{v,u,w\}.
\end{eqnarray*}
Furthermore, $[\c,\c,\c]$ is regarded as
\begin{eqnarray*}
 [u,v,w]_{C}=\rho(Tu,Tv)w+\rho(Tv,Tw)u+\rho(Tw,Tu)v.
\end{eqnarray*}
Since $T$ is an $\mathcal{O}$-operator, we have
\begin{eqnarray*}
T[u,v,w]_C=[Tu,Tv,Tw].
\end{eqnarray*}
For any $v_1,v_2, v_3, v_4, v_5\in V$, we have
\begin{eqnarray*}
\{\b(v_1), \b(v_2), \{v_3, v_4, v_5\}\}&=& \rho(T\b(v_1), T\b(v_2)) \rho(Tv_3, Tv_4)v_5,\\
\{[v_1,v_2,v_3]_L, \b(v_4), \b(v_5)\}&=& \rho(T[v_1,v_2,v_3]_L, T\b(v_4))\b(v_5)=\rho([Tv_1,Tv_2,Tv_3], T\b(v_4))\b(v_5),\\
\{\b(v_3), [v_1,v_2,v_4]_L, \b(v_5)\}&=& \rho(T\b(v_3), T[v_1,v_2,v_4]_L)\b(v_5)=\rho(T\b(v_3), [Tv_1,Tv_2,Tv_4])\b(v_5),\\
\{\b(v_3), \b(v_4), \{v_1, v_2, v_5\}\}&=& \rho(T\b(v_3), T\b(v_4)) \rho(Tv_1, Tv_2)v_5.
\end{eqnarray*}
Since $(V, \rho, \b)$ is a representation, we can check that (3.4) and (3.5) hold. And this finishes the proof.

\begin{proposition}
Let $(L, [\c, \c, \c]_L, \a)$ be a 3-Hom-Lie algebra.
Then there exists a  compatible 3-Hom-pre-Lie algebra if and only if there exists an invertible $\mathcal{O}$-operator on $L$.
\end{proposition}
{\bf Proof.} Let $T$ be an invertible $\mathcal{O}$-operator of $L$ associated to a representation $(V, \rho, \b)$.
 Then there exists a 3-Hom-pre-Lie algebra structure on $V$ defined by
\begin{eqnarray*}
\{u,v,w\}=\rho(Tv,Tw)u, ~~~\forall u,v,w\in V.
\end{eqnarray*}
Moreover, there is an induced 3-Hom-pre-Lie algebra structure $\{\c, \c, \c\}_L$ on $L=T(V)$ given by
\begin{eqnarray*}
\{x, y, z\}_L=T\{T^{-1}x, T^{-1}y, T^{-1}z\}=T\rho(x,y)T^{-1}z.
\end{eqnarray*}
Since $T$ is an $\mathcal{O}$-operator, we have
\begin{eqnarray*}
[x,y, z]&=& T\rho(y,z)T^{-1}x+T\rho(z,x)T^{-1}y+T\rho(x,y)T^{-1}z\\
&=& \{x, y, z\}_L+\{y, z, z\}_L+\{z, x, y\}_L.
\end{eqnarray*}
Therefore, $(L, \{\c, \c, \c\}_L, \a)$ is a compatible 3-Hom-pre-Lie algebra. Conversely, the identity map $id$ is an $\mathcal{O}$-operator on $L$.\hfill $\square$
\begin{definition}(\cite{WW17})
Let $(L, [\c, \c, \c], \a)$ be 3-Hom-Lie algebra and $r\in L\o L$. The equation
\begin{eqnarray*}
[[r,r,r]]=0
\end{eqnarray*}
is called the 3-Lie classical Hom-Yang-Baxter equation.
\end{definition}

For any $r\in L\o L$, the induced skew-symmetric linear map $r: L^{\ast}\rightarrow L$ is defined by
\begin{eqnarray*}
<r(\xi), \eta>=<r, \xi\eta>.
\end{eqnarray*}
We denote the ternary operation $\Delta^{\ast}: L^{\ast}\o L^{\ast}\o  L^{\ast}\rightarrow L^{\ast}$ by $[\c,\c,\c]^{\ast}$.
According to \cite{WW17}, for any $r=\sum_i x_i\o y_i\in L\o L$ and $x\in L$, we have
\begin{eqnarray*}
&& \Delta_1(x)=\sum_{i,j}[x,x_i,x_j]\o \a(y_j)\o \a(y_i),\\
&& \Delta_2(x)=\sum_{i,j}\a(y_i)\o [x,x_i,x_j]\o \a(y_j),\\
&& \Delta_3(x)=\sum_{i,j}\a(y_j) \o \a(y_i) \o [x,x_i,x_j].
\end{eqnarray*}

\begin{proposition}
Let $(L, [\c, \c, \c], \a)$ be a 3-Hom-Lie algebra and $r\in L\o L$ such that $\a^{\otimes^2}r=r$.
Suppose that $r$ is skew-symmetric and $\Delta=\Delta_1+\Delta_2+\Delta_3: L\rightarrow L\o L\o L$.
Then we have
\begin{eqnarray}
[\xi,\eta,\gamma]^{\ast}=ad^{\ast}_{r(\xi),r(\eta)}\g+ad^{\ast}_{r(\eta),r(\g)}\xi+ad^{\ast}_{r(\g),r(\xi)}\eta.
\end{eqnarray}
Furthermore, we have
\begin{eqnarray}
[r(\xi),r(\eta),r(\g)]-r([\xi,\eta,\gamma]^{\ast})=[[r,r,r]](\xi,\eta,\gamma),
\end{eqnarray}
for any $\xi,\eta,\gamma\in L^{\ast}$.
\end{proposition}
{\bf Proof.}    Let $r=\sum_{i}x_i\o y_i$, for any $x,y\in L$ and $\xi, \eta, \gamma \in L^{\ast}$, we have
\begin{eqnarray*}
<x, ad^{\ast}_{r\circ \a^{\ast}(\xi),r\circ \a^{\ast}(\eta) }\gamma>&=& <-[r\circ \a^{\ast}(\xi),r\circ \a^{\ast}(\eta), x],\g>\\
&=& -<r, \a^{\ast}(\eta)\o ad^{\ast}_{r\circ \a^{\ast}(\xi),x}\g>\\
&=& \sum_i <y_i, \a^{\ast}(\eta)> <r, \a^{\ast}(\xi) \o ad^{\ast}_{x,x_i}\g>\\
&=& \sum_i <y_i, \a^{\ast}(\eta)><y_j, \a^{\ast}(\xi)><[x,x_i,x_j],\g>\\
&=& \sum_{i,j} <\a(y_i), \eta><\a(y_j), \xi><[x,x_i,x_j],\g>\\
&=&<\sum_{i,j} \a(y_j)\o \a(y_i)\o [x,x_i,x_j], \xi\o \eta\o \g>\\
&=& <\Delta_3(x), \xi\o \eta\o \g>.
\end{eqnarray*}
Similarly, we have
\begin{eqnarray*}
<x, ad^{\ast}_{r\circ \a^{\ast}(\eta),r\circ \a^{\ast}(\gamma) }\xi>= <\Delta_1(x), \xi\o \eta\o \g>,
<x, ad^{\ast}_{r\circ \a^{\ast}(\g),r\circ \a^{\ast}(\xi) }\eta>= <\Delta_1(x), \xi\o \eta\o \g>.
\end{eqnarray*}
Therefore, we have
\begin{eqnarray*}
&&<\Delta(x), \xi\o \eta\o \g>\\
&=& <\Delta_1(x)+\Delta_2(x)+\Delta_3(x), \xi\o \eta\o \g>\\
&=& <x, ad^{\ast}_{r\circ \a^{\ast}(\eta),r\circ \a^{\ast}(\gamma) }\xi>+<x, ad^{\ast}_{r\circ \a^{\ast}(\g),r\circ \a^{\ast}(\xi) }\eta>+<x, ad^{\ast}_{r\circ \a^{\ast}(\xi),r\circ \a^{\ast}(\eta) }\gamma>\\
&=& <x, [\xi,\eta, \gamma]^{\ast}>.
\end{eqnarray*}
So Eq. (3.7) holds as required.
For  Eq. (3.8)  we take any $\kappa\in L^{\ast}$ and compute
\begin{eqnarray*}
&&[[r,r,r]](\xi,\eta,\gamma, \kappa)\\
&=& \sum_{i,j,k} ([x_i,x_j,x_k]\o \a(y_i)\o \a(y_j)\o \a(y_k)(\xi,\eta,\gamma, \kappa)+\a(x_i)\o [y_i,x_j,x_k]\o \a(y_j)\o \a(y_k)(\xi,\eta,\gamma, \kappa)\\
&& \a(x_i)\o \a(x_j)\o [y_i,y_j,x_k]\o \a(y_k)(\xi,\eta,\gamma, \kappa)+\a(x_i)\o \a(x_j)\o \a(x_k)\o [y_i,y_j,y_k](\xi,\eta,\gamma, \kappa))\\
&=& \sum_{i,j,k}<\xi, [x_i,x_j,x_k]><\eta, \a(y_i)> <\gamma, \a(y_j)> <\kappa, \a(y_k)>+<\eta, [y_i,x_j,x_k]> <\xi, \a(x_i)>\\
&& <\g, \a(y_j)> <\kappa, \a(y_k)> <\g, [y_i,y_j,x_k]><\xi, \a(x_i)> <\eta, \a(x_j)> <\kappa, \a(y_k)>+\\
&&<\kappa, [y_i,y_j,y_k]> <\xi, \a(x_i)> <\eta, \a(x_j)> <\g, \a(x_k)>
\end{eqnarray*}
\begin{eqnarray*}
&=&-<\xi, [r\circ \a^{\ast}(\eta),r\circ \a^{\ast}(\gamma), r\circ \a^{\ast}(\kappa) >-<\eta, [r\circ \a^{\ast}(\g),r\circ \a^{\ast}(\xi), r\circ \a^{\ast}(\kappa) >\\
&&-<\g, [r\circ \a^{\ast}(\xi),r\circ \a^{\ast}(\eta), r\circ \a^{\ast}(\kappa) >+<\kappa, [r\circ \a^{\ast}(\xi),r\circ \a^{\ast}(\eta), r\circ \a^{\ast}(\gamma) >\\
&=& <[r\circ \a^{\ast}(\xi),r\circ \a^{\ast}(\eta), r\circ \a^{\ast}(\gamma)]-r\circ \a^{\ast}([\xi, \eta, \gamma]^{\ast}), \kappa>
\end{eqnarray*}
So Eq. (3.8) holds and this finishes the proof.
\hfill $\square$

\begin{proposition}
Let $(L, [\c, \c, \c]_L, \a)$ be a regular 3-Hom-Lie algebra and $r\in L\o L$ such that $\a^{\otimes^2}r=r$. Suppose $r$ is skew-symmetric and nondegenerate. Then $r$ is a solution of the 3-Hom-Lie Yang-Baxter equation if and only if nondegenerate skew-symmetric bilinear form $B$ on $L$ defined by $B(x,y)=<r^{-1}(x), y>$ satisfies
\begin{eqnarray}
B(\a[x,y,z],w)- B(\a[x,y,w],z)+B(\a[x,z,w],y)-B(\a[y,z,w],x)=0,
\end{eqnarray}
for any $x,y,z, w\in L.$
\end{proposition}
{\bf Proof.}  For any $x,y,z, w\in L$, there exists $\xi, \eta, \gamma, \kappa \in L^{\ast}$ such that  $r(\xi)=x, r(\eta)=y, r(\gamma)=z, r(\kappa)=w$.
If $[[r,r,r]]=0$, we have
\begin{eqnarray*}
&&B(\a[x,y,z],w)\\
&=& <\a[r(\xi), r(\eta),r(\gamma)], \kappa>\\
&=& <r\circ \a^{\ast}(ad^{\ast}_{r\circ \a^{\ast}(\xi), r\circ \a^{\ast}(\eta)}\g+ad^{\ast}_{r\circ \a^{\ast}(\eta), r\circ \a^{\ast}(\gamma)}\xi+ad^{\ast}_{r\circ \a^{\ast}(\gamma), r\circ \a^{\ast}(\xi)}\eta), \kappa>\\
&=&<ad^{\ast}_{r\circ \a^{\ast}(\xi), r\circ \a^{\ast}(\eta)}\g+ad^{\ast}_{r\circ \a^{\ast}(\eta), r\circ \a^{\ast}(\gamma)}\xi+ad^{\ast}_{r\circ \a^{\ast}(\gamma), r\circ \a^{\ast}(\xi)}\eta, -\a \circ r(\kappa)>\\
&=& <\g,\a[x,y,w]>-<-\xi,\a[y,z,w]>-<-\eta,\a[z,x,w]>\\
&=& B(\a[x,y,w],z)-B(\a[x,z,w],y)+B(\a[y,z,w],x).
\end{eqnarray*}
Thus the proof is finished. \hfill $\square$

\section{Representations of 3-Hom-pre-Lie algebras}
\def\theequation{\arabic{section}.\arabic{equation}}
\setcounter{equation} {0}

In this section, we introduce the notion of the representation of a 3-Hom-pre-Lie algebra and construct the
corresponding semidirect product 3-Hom-pre-Lie algebra.

\begin{definition}
A representation of a 3-Hom-pre-Lie algebra $(L, \{\c, \c, \c\}, \a)$ on a vector space $V$ consists
of a pair $(\rho, \mu)$, where $\rho :  L\wedge L\rightarrow gl(V )$ is a representation of the 3-Hom-Lie algebra $(L, [\c, \c, \c]_{L}, \a)$ on $V$ and
$\mu: L\o L\rightarrow gl(V )$ is a linear map such that the following equalities hold:
\begin{eqnarray}
\rho(x_1,x_2)\mu(x_3,x_4)&=&\mu(x_3,x_4) \rho(x_1,x_2)-\mu(x_3,x_4) \mu(x_2,x_1)\nonumber\\
&&\mu(x_3,x_4) \mu(x_1,x_2) +\mu([x_1,x_2,x_3]_L, x_4)+\mu(x_3, \{x_1, x_2,x_4\}),~~~~~~\\
\mu([x_1,x_2,x_3]_C, x_4)&=&\rho(x_1,x_2)\mu(x_3,x_4)+\rho(x_2,x_3)\mu(x_1,x_4)+\rho(x_3,x_1)\mu(x_2,x_4),~~~~~\\
 \mu(x_1, \{x_1,x_3,x_4\})&=&\mu(x_3,x_4) \mu(x_1,x_2)+\mu(x_3,x_4) \rho(x_1,x_2)\nonumber\\
&&-\mu(x_3,x_4) \mu(x_2,x_1)-\mu(x_2,x_4) \mu(x_1,x_3)\nonumber\\
&& -\mu(x_2,x_4)\rho(x_1,x_3)+\mu(x_2,x_4) \mu(x_3,x_1)+\rho(x_2,x_3) \mu(x_1,x_4),
\end{eqnarray}
\begin{eqnarray}
\mu(x_3,x_4) \rho(x_1,x_2)&=&\mu(x_3,x_4) \mu(x_2,x_1)-\mu(x_3,x_4) \mu(x_1,x_2)\nonumber\\
&& +\rho(x_1,x_2) \rho(x_3,x_4)-\mu(x_2, \{x_1,x_3,x_4\})+\mu(x_1, \{x_2,x_3,x_4\}),~~~~~
\end{eqnarray}
for all $x_1,x_2,x_3,x_4\in L$.
\end{definition}

\begin{theorem}
Let $(V, \rho, \mu)$ be a representation of a 3-Hom-pre-Lie algebra $(L, \{\c, \c, \c\}, \a)$. Define a trilinear bracket
operation $\{\c,\c,\c\}_{\rho,\mu}: (L\oplus V)\o (L\oplus V)\o (L\oplus V)\rightarrow (L\oplus V)$ by
\begin{eqnarray*}
\{x_1+v_1,x_2+v_2,x_3+v_3\}_{\rho,\mu}=\{x_1,x_2,x_3\}+\rho(x_1,x_2)v_3+\mu(x_2,x_3)v_1-\mu(x_1,x_3)v_2,\\
(\a+\b)(x_1+v_1)=\a(x_1)+\b(v_1),
\end{eqnarray*}
for any $x_1,x_2,x_3\in L$ and  $v_1,v_2,v_3\in V$.  Then $(L\oplus V, \{\c,\c,\c\}_{\rho,\mu}, \a+\b)$ is a 3-Hom-pre-Lie algebra.
\end{theorem}

{\bf Proof.}   Straightforward. \hfill $\square$

\medskip

Let $V$ be a vector space. Define the switching operator $\tau: V\o V\rightarrow V\o V$ by
\begin{eqnarray*}
\tau(v_1\o v_2)=v_2\o v_1,~~~\forall v_1, v_2\in V.
\end{eqnarray*}
\begin{proposition}
Let $(\rho, \mu)$ be a representation of a 3-Hom-pre-Lie algebra $(L, \{\c, \c, \c\}, \a)$ on a vector space $V$. Then $\rho-\mu\tau+\mu$ is a representation of the sub-adjacent 3-Hom-Lie algebra $(L^c, [\c, \c, \c]_{C}, \a)$ on the vector space $V$.
\end{proposition}
{\bf Proof.} Consider the sub-adjacent 3-Hom-Lie algebra structure $[\c,\c,\c]_L$, for any $x_1,x_2,x_3$ and $v_1,v_2,v_3\in V$, we have
\begin{eqnarray*}
&&[x_1+v_1,x_2+v_2,x_3+v_3]_L\\
&=& \{x_1+v_1,x_2+v_2,x_3+v_3\}_{\rho,\mu}+\{x_2+v_2,x_3+v_3, x_1+v_1\}_{\rho,\mu}\\
&&+\{x_3+v_3, x_1+v_1,x_2+v_2\}_{\rho,\mu}\\
&=& \{x_1,x_2,x_3\}+\rho(x_1,x_2)v_3+\mu(x_2,x_3)v_1-\mu(x_1,x_3)v_2\\
&&\{x_2,x_3,x_1\}+\rho(x_2,x_3)v_1+\mu(x_3,x_1)v_2-\mu(x_2,x_1)v_3\\
&&\{x_3,x_1,x_2\}+\rho(x_3,x_1)v_2+\mu(x_1,x_2)v_3-\mu(x_3,x_2)v_1\\
&=& [x_1,x_2,x_3]_L+((\rho-\mu\tau+\mu)(x_1,x_2))v_3\\
&& ((\rho-\mu\tau+\mu)(x_2,x_3))v_1+((\rho-\mu\tau+\mu)(x_3,x_1))v_2.
\end{eqnarray*}
By Lemma 1.3, $\rho-\mu\tau+\mu$ is a representation of the sub-adjacent 3-Hom-Lie algebra $(L^c, [\c, \c, \c]_{C}, \a)$ on the vector space $(V, A)$.  And the proof is finished. \hfill $\square$

\begin{corollary}
Let $(\rho, \mu)$ be a representation of a 3-Hom-pre-Lie algebra $(L, \{\c,\c,\c\}, \a)$ on a vector space
$V$. Then the semidirect product 3-Hom-pre-Lie algebras $L\ltimes_{\rho, \mu}V$ and $L\ltimes_{\rho-\mu\tau+\mu, 0}V$ given by the
representations $(\rho, \mu)$ and $(\rho-\mu\tau+\mu, 0)$ respectively have the same sub-adjacent 3-Hom-Lie algebra
$L^c\ltimes_{\rho-\mu\tau+\mu}V$, which is the semidirect product of the 3-Hom-Lie algebra $(L^c, [\c, \c, \c]_{C}, \a)$ and
its representation $(V, \rho-\mu\tau+\mu)$.
\end{corollary}
\begin{proposition}
Let $(\rho, \mu)$ be a representation of a 3-Hom-pre-Lie algebra $(L, \{\c, \c, \c\}, \a)$ on a vector space $(V, A)$. Then $(\rho^{\ast}-\mu^{\ast}\tau+\mu^{\ast}, -\mu^{\ast})$ is a representation of 3-Hom-pre-Lie algebra $(L, \{\c, \c, \c\}, \a)$ on a vector space $(V^{\ast}, A^{\ast})$.
\end{proposition}
{\bf Proof.}  By Proposition 4.3, $\rho-\mu\tau+\mu$ is a representation of the sub-adjacent 3-Hom-Lie algebra $(L, [\c,\c,\c]_L, \a)$ on the vector space $(V, A)$. It is easy to check that $\rho^{\ast}-\mu^{\ast}\tau+\mu^{\ast}$ is a representation of the sub-adjacent 3-Hom-Lie algebra $(L^c, [\c, \c, \c]_{C}, \a)$ on the dual vector space $(V^{\ast}, A^{\ast})$.  And the proof is finished. \hfill $\square$

\begin{corollary}
Let $(V, \rho, \mu)$ be a representation of a 3-Hom-pre-Lie algebra $(L, \{\c, \c, \c\}, \a)$. Then the
semidirect product 3-Hom-pre-Lie algebras $L\ltimes_{\rho^{\ast}, 0}V^{\ast}$ and $L\ltimes_{\rho^{\ast}-\mu^{\ast}\tau+\mu^{\ast}, -\mu^{\ast}}V^{\ast}$  given by the representations $(\rho^{\ast}, 0)$ and $(\rho^{\ast}-\mu^{\ast}\tau+\mu^{\ast}, -\mu^{\ast})$ respectively have the same sub-adjacent 3-Hom-Lie algebra $L\ltimes_{\rho^{\ast}}V^{\ast}$, which is the semidirect product of the 3-Hom-Lie algebra $(L^c, [\c, \c, \c]_{C}, \a)$ and its representation $(V, \rho^{\ast})$.
\end{corollary}
\section{Symplectic structures and phase spaces of 3-Hom-Lie algebras}
\def\theequation{\arabic{section}.\arabic{equation}}
\setcounter{equation} {0}

\begin{definition}
A symplectic structure on a regular 3-Hom-Lie algebra $(L, [\c, \c, \c], \a)$ is a nondegenerate
skew-symmetric bilinear form $\omega\in L^{\ast}\wedge L^{\ast}$ such that $\omega\circ \a=\omega$ satisfying the following equality
\begin{eqnarray}
\omega([x,y,z], \a(w))-\omega([y,z,w], \a(x))+\omega([z,w,x], \a(y))-\omega([w,x,y], \a(z))=0,
\end{eqnarray}
for any $ x,y,z,w\in L$.
\end{definition}

\begin{definition}$^{\cite{SC15}}$
Let $(L, [\c, \c, \c], \a)$ be a 3-Hom-Lie algebra. $B: L\times L\rightarrow F$ be a non-degenerate symmetric bilinear form on $L$. If $B$ satisfies
\begin{eqnarray}
B([x,y,z],w)+B(z, [x,y,w])=0,~~\forall x,y,z,w\in L.
\end{eqnarray}
Then $B$ is called a a metric on 3-Hom-Lie algebra $(L, [\c, \c, \c], \a)$, and $(L, [\c, \c, \c], \a, B)$ is a metric 3-Hom-Lie algebra.

\end{definition}

If there exists a metric $B$ and a symplectic structure $\omega$ on 3-Hom-Lie algebra $(L, [\c, \c, \c], \a)$, then $(L, [\c, \c, \c], \a, B, \omega)$ is called a metric symplectic 3-Hom-Lie algebra.
\medskip

Let $(L, [\c, \c, \c], \a, B)$ be a metric 3-Hom-Lie algebra, we denote
\begin{eqnarray*}
Der_B(L)=\{D\in Der(L)|B(Dx,y)+B(x,Dy)=0, \forall x,y\in L\}=Der(L) \cap so(L, B).
\end{eqnarray*}
\begin{theorem}
Let $(L, [\c, \c, \c], \a, B)$ be a metric 3-Hom-Lie algebra. Then there exists a symplectic structure on $L$ if and only if there exists a skew-symmetric invertible derivation $D\in Der_B(L)$.

\end{theorem}
{\bf Proof.} For any  $x,y\in L$, define $D:L \rightarrow L$ by
\begin{eqnarray}
B(Dx,y)=\omega (\a(x),y).
\end{eqnarray}
It is clear that $D$ is invertible. Next we will check that $D$ is a skew-symmetric invertible derivation of $(L, [\c, \c, \c], \a, B)$. In fact, for any $x,y,z,w\in L$, we have
\begin{eqnarray*}
&&B([Dx,y,z], w)+B([x,Dy,z],w)+ B([x,y,Dz], w)+B(D[x,y,z],w)\\
&=& -B([y,z,w], Dx)+B([x,z,w], Dy)-B([x, y,w], Dz)+B([x, y,z], Dw)\\
&=& \omega([x,y,z], \a(w))-\omega([y,z,w], \a(x))+\omega([z,w,x], \a(y))-\omega([w,x,y], \a(z))=0,
\end{eqnarray*}
that is, $D\in Der_B(L)$.

Conversely, assume $D\in Der_B(L)$. Define $\omega$ by (5.3), then there exists a symplectic structure on $L$ satisfies (5.1).
\begin{example}
Let $(L, [\c, \c, \c], \a)$ be a 3-Hom-Lie algebra and
\begin{eqnarray*}
F[t]=\{f(t)=\sum^{m}_{i=0}a_it^{i}|a_i\in F,m\in N\}
\end{eqnarray*}
the algebra of polynomials over $F$. We consider
\begin{eqnarray*}
L_n=L\o (tF[t]/t^{n}F[t]),
\end{eqnarray*}
where $tF[t]/t^{n}F[t]$ is the quotient space of $tF[t]$ module $t^nF[t]$. Then $L_n$ is a nilpontent 3-Hom-Lie algebra in the following multiplication
\begin{eqnarray*}
\a(x\o t^{\overline{p}})=\a(x)\o t^{\overline{p}},
[x\o t^{\overline{p}},y\o t^{\overline{q}}, z\o t^{\overline{r}} ]=[x,y,z]_{L}\o t^{\overline{p+q+r}},
\end{eqnarray*}
for any $x,y,z\in L,p,q,r\in N\ \{0\}.$
Define an endomorphism $D$ of $L_n$ by
\begin{eqnarray*}
D(x\o t^{\overline{p}})=p(x\o t^{\overline{p}}), \forall x\in L, p=1,...,n-1.
\end{eqnarray*}
Then $D$ ia an invertible derivation of the 3-Hom-Lie algebra $L_n$.

Let $\widetilde{L}_{n}=L_n\oplus L^{\ast}_n$, where $L^{\ast}_n$ is the dual space of $L_n$. Then $(\widetilde{L}_{n}, B)$ ia a metric 3-Hom-Lie algebra with the multiplication
\begin{eqnarray*}
&&[x+f, y+g,z+h]=[x,y,z]+ad^{\ast}(y,z)f-ad^{\ast}(x,z)g+ad^{\ast}(x,y)h, \\
&& B(x+f, y+g)=f(y)+g(x),
\end{eqnarray*}
for any $x,y,z\in L_n, f,g,h\in L^{\ast}_n$.
And define linear maps $\widehat{D}, \widetilde{\a}: \widetilde{L}_{n}\rightarrow \widetilde{L}_{n}$ by
\begin{eqnarray*}
\widehat{D}(x+f)=Dx+D^{\ast}f,
\widetilde{\a}(x+f)=\a(x)+f\circ \a,
\end{eqnarray*}
where $D^{\ast}f=-fD$. Then $\widehat{D}$ is an invertible. Hence $(\widetilde{L}_{n}, \widetilde{\a}, B, \omega)$ is a metric symplectic 3-Hom-Lie algebra,where $\omega$ is defined as follows:
\begin{eqnarray*}
\omega(\widetilde{\a}(x+f),y+g)=B(\widehat{D}(x+f), y+g)=-f(Dy)+g(Dx).
\end{eqnarray*}
\end{example}

\begin{proposition}
Let $(L, [\c, \c, \c], \a, \omega)$ be a symplectic 3-Hom-Lie algebra. Then there exists a compatible 3-Hom-pre-Lie algebra structure $\{\c,\c,\c\}$ on $L$ given by
\begin{eqnarray}
\omega(\{x,y,z\}, \a(w))=-\omega(\a(z), [x,y,w]), \forall x,y,z,w\in L.
\end{eqnarray}
\end{proposition}

{\bf Proof.}  For any $x,y\in L$, define the map $T: L^{\ast}\rightarrow L$ by $<T^{-1}x, y>=\omega(x,y)$.
 By (5.1), we obtain that $T$ is an invertible $\mathcal{O}$-operator associated to the coadjoint representation $(L^{\ast}, ad^{\ast}, \a^{\ast})$,  there exists a compatible  3-Hom-pre-Lie algebra on $L$ given by $\{x, y, z\}=T(ad^{\ast}_{x,y}T^{-1}z)$, for any $x,y,z\in L$, we have
\begin{eqnarray*}
\omega(\{x,y,z\}, \a(w))&=&\omega(T(ad^{\ast}_{x,y}T^{-1}z), \a(w))=<ad^{\ast}_{x,y}T^{-1}z, \a(w)>\\
&=& <T^{-1}(\a(z)), -[x,y,w]>=-\omega(\a(z), [x,y,w]).
\end{eqnarray*}
The proof is finished.  \hfill $\square$
\medskip

Let $V$ be a vector space and $V^{\ast}$ its dual space. Then there is a natural nondegenerate
skew-symmetric bilinear form $\omega$ on $T ^{\ast}V = V \oplus V^{\ast}$ given by:
\begin{eqnarray}
\omega(x+f, y+g) =<f, y>-<g, x>, \forall x,y\in V, f,g\in V^{\ast}.
\end{eqnarray}

\begin{definition}
Let $(L, [\c, \c, \c], \a)$ be a  3-Hom-Lie algebra and $L^{\ast}$ its dual space.
If there is a 3-Hom-Lie algebra structure $[\c, \c, \c]$ on the direct sum vector space $T ^{\ast}L = L \oplus L^{\ast}$ such
that $(L \oplus L^{\ast}, [\c,\c, \c], \a+\a^{\ast}, \omega)$ is a symplectic 3-Hom-Lie algebra, where $\omega$ given by (5.5),
and $(L, [\c, \c, \c], \a)$
and $(L^{\ast}, [\c, \c, \c]^{\ast}, \a^{\ast})$ are 3-Hom-Lie subalgebras of $(L \oplus L^{\ast}, [\c,\c, \c], \a+\a^{\ast})$,
then the symplectic 3-Hom-Lie algebra $(L \oplus L^{\ast}, [\c,\c, \c], \a+\a^{\ast}, \omega)$ is called a phase space of the 3-Hom-Lie algebra $(L, [\c, \c, \c], \a)$.
\end{definition}

Next, we study the relation between phase spaces of 3-Hom-Lie algebra and 3-Hom-pre-Lie algebra.
\begin{theorem}
A 3-Hom-Lie algebra has a phase space if and only if it is sub-adjacent to a 3-Hom-pre-Lie
algebra.
\end{theorem}
{\bf Proof.} $\Leftarrow$  Assume $(L, [\c, \c, \c], \a)$ is a  3-Hom-pre-Lie algebra. By Lemma 1.3 and  Proposition 3.4,    the left multiplication $\mathcal{L}$ is a
representation of the sub-adjacent 3-Lie algebra $L^{C}$ on $L$,  $\mathcal{L^{\ast}}$  is a representation
of the sub-adjacent 3-Lie algebra $L^{C}$ on $L^{\ast}$,  then  we have a 3-Hom-Lie algebra $(L^{C}\oplus L^{\ast}, [\c,\c,\c]_{L^{\ast}}, \a+\a^{\ast})$.
For all $x_1,x_2,x_3,x_4\in L$ and $f_1,f_2,f_3,f_4\in L^{\ast}$, we have
\begin{eqnarray*}
&&\omega([x_1+f_1,x_2+f_2,x_3+f_3]_{L^{\ast}}, \a(x_4)+\a^{\ast}(f_4))\\
&=&\omega([x_1,x_2,x_3]_C+ \mathcal{L^{\ast}}(x_1,x_3)f_3+\mathcal{L^{\ast}}(x_2,x_3)f_1 +\mathcal{L^{\ast}}(x_3,x_1)f_2, \a(x_4)+\a^{\ast}(f_4))\\
&=& <\mathcal{L^{\ast}}(x_1,x_3)f_3+\mathcal{L^{\ast}}(x_2,x_3)f_1 +\mathcal{L^{\ast}}(x_3,x_1)f_2,\a(x_4)>-<\a^{\ast}(f_4), [x_1,x_2,x_3]_C> \\
&=& -<\a^{\ast}(f_3), \{x_1,x_2,x_3\}>-<\a^{\ast}(f_1), \{x_2,x_3,x_4\}>-<\a^{\ast}(f_2), \{x_3,x_1,x_4\}>\\
&& -<\a^{\ast}(f_4), \{x_1,x_2,x_3\}>-<\a^{\ast}(f_4), \{x_2,x_3,x_1\}>-<\a^{\ast}(f_4), \{x_3,x_1,x_2\}>.
\end{eqnarray*}
Similarly, we have
\begin{eqnarray*}
&&\omega([x_2+f_2,x_3+f_3,x_4+f_4]_{L^{\ast}}, \a(x_1)+\a^{\ast}(f_1))\\
&=& -<\a^{\ast}(f_4), \{x_2,x_3,x_1\}>-<\a^{\ast}(f_2), \{x_3,x_4,x_1\}>-<\a^{\ast}(f_3), \{x_4,x_2,x_1\}>\\
&& -<\a^{\ast}(f_1), \{x_2,x_3,x_4\}>-<\a^{\ast}(f_1), \{x_3,x_4,x_2\}>-<\a^{\ast}(f_1), \{x_4,x_2,x_3\}>,\\
&&\omega([x_3+f_3,x_4+f_4, x_1+f_1,]_{L^{\ast}}, \a(x_2)+\a^{\ast}(f_2))\\
&=& -<\a^{\ast}(f_1), \{x_3,x_4,x_2\}>-<\a^{\ast}(f_3), \{x_4,x_1,x_2\}>-<\a^{\ast}(f_4), \{x_1,x_3,x_2\}>\\
&& -<\a^{\ast}(f_2), \{x_3,x_4, x_1\}>-<\a^{\ast}(f_2), \{x_4,x_1,x_3\}>-<\a^{\ast}(f_2), \{x_1,x_3,x_4\}>,\\
&&\omega([x_4+f_4, x_1+f_1,x_2+f_3,]_{L^{\ast}}, \a(x_3)+\a^{\ast}(f_3))\\
&=& -<\a^{\ast}(f_2), \{x_4,x_1,x_3\}>-<\a^{\ast}(f_4), \{x_1,x_2,x_3\}>-<\a^{\ast}(f_1), \{x_2,x_4,x_3\}>\\
&& -<\a^{\ast}(f_3), \{x_4,x_1,x_2\}>-<\a^{\ast}(f_3), \{x_1,x_2,x_4\}>-<\a^{\ast}(f_3), \{x_2,x_4,x_1\}>,
\end{eqnarray*}
So $\omega$ is a symplectic structure on the semidirect product 3-Hom-Lie algebra   $(L^{C}\oplus L^{\ast},$\\
$ [\c,\c,\c]_{L^{\ast}}, \a+\a^{\ast})$. Thus the symplectic 3-Lie algebra $(L^{C}\oplus L^{\ast}, [\c,\c,\c]_{L^{\ast}}, \a+\a^{\ast}, \omega)$ is a phase space of the sub-adjacent 3-Hom-Lie algebra $(L^{C}, [\c,\c,\c]_C, \a)$.

$\Rightarrow$  Clearly. \hfill $\square$
\begin{center}
 {\bf ACKNOWLEDGEMENT}
\end{center}

  The work of S. J. Guo is  supported by the NSF of China (No. 11761017) and the Youth Project for Natural Science Foundation of Guizhou provincial department of education (No. KY[2018]155).
  The work of X. H. Zhang is supported by the NSF of China (No. 11801304) and the Project Funded by China Postdoctoral Science Foundation (No. 2018M630768).
   The work of S. X. Wang is  supported by  the outstanding top-notch talent cultivation project of Anhui Province (No. gxfx2017123)
 and the Anhui Provincial Natural Science Foundation (No. 1808085MA14).

\end{document}